\input amstex
\documentstyle{amsppt}
\TagsOnRight
\pagewidth{5.5 in}
\pageheight{8 in}
\magnification=1200
\NoBlackBoxes
 

\def\inner#1#2#3{ \langle #1,\, #2 \rangle_{_{#3}}}

\def\sr{\Cal S(\Bbb R)}
\def\vartau{\tau\hskip2pt\!\!\!\!\iota\hskip1pt}
\def\ft#1{\widehat{#1}}

\def\intg{\int_{-\infty}^\infty}
 
\def\BW{1}
\def\BEEKa{2}
\def\BEEKb{3}
\def\BK{4}
\def\AC{5}
\def\MRa{6}
\def\MRb{7}
\def\WR{8}
\def\SWa{9}
\def\SWb{10}

%
%
 
\topmatter
\pretitle{\hfill \sevenit C.~R.~Math.~Rep.~Acad.~Sci.~Canada (to appear)}

\vskip20pt
 
\title
Periodic Integral Transforms \\ and C*-algebras
\endtitle
 
\rightheadtext{Periodic Integral Transforms}
 
\author
S.~Walters
\endauthor

\address
Department of Mathematics, 
University of Northern British Columbia,
Prince George, B.C.  V2N 4Z9  CANADA
\endaddress
 
\email
walters\@hilbert.unbc.ca \ or \ walters\@unbc.ca  \hfill
\break\indent{\it Home page}: http://hilbert.unbc.ca/walters
\endemail
 
\thanks
Research partly supported by NSERC grant OGP0169928 \hfill {\sevenrm (\TeX
File: transform36.tex)}
\endthanks
 
 

\abstract
We construct canonical integral transforms, analogous to the Fourier transform,
that have periods six and three.  The existence of this transform is shown to 
arise naturally from the expectation that the Schwartz space on the real line,
viewed as the Heisenberg module of Rieffel and Connes over the rotation 
C*-algebra, should extend to a module action over the crossed product of the
latter by the canonical automorphisms of orders three and six (which does in
fact happen and is shown here).
\endabstract

\endtopmatter


\document

\subhead \S1. INTRODUCTION \endsubhead

It is a well known classical fact that the Fourier transform of
a Schwartz function $f$
$$
\ft f(t) = \intg f(x) e(-tx) dx, \tag1
$$
has period four and extends to a unitary operator on $L^2(\Bbb R)$.
(Throughout the paper we write $e(t) := e^{2\pi it}$.)
This stems from the fact that $\Hat{\Hat f} (t) = f(-t)$.  In this paper
we show that if the product $tx$ in (1) is replaced by a suitable quadratic,
then one obtains transforms of period three and six.  More specifically,
one has a one-parameter family of {\it hexic} transforms
$$
(Hf)(t) = i^{1/6}\sqrt{2\mu} \intg f(x) e(2\mu tx-\mu x^2) dx, \tag2
$$
for $\mu>0$ and $f$ in the Schwartz space $\sr$.  Thus, $H$ extends to a unitary 
operator on $L^2(\Bbb R)$ of period six (i.e., $H^6=I$).
(The ``ideal'' transform is when $\mu=\tfrac12$ as is explained in Remark 2 below.)
Note that $H$ is a composition of the multiplication operator by the complex Gaussian 
$e(-\mu x^2)$ and an inverse Fourier transform (up to scaling), and hence is itself 
a unitary operator on $L^2(\Bbb R)$ that leaves invariant $\sr$.  

\proclaim{Theorem 1}
One has $(H^3f)(t)=f(-t)$ for all $f \in \sr$, so that the transform $H$ has 
period six and extends to a unitary operator on $L^2(\Bbb R)$.
Further, its square $H^2$ (the cubic transform) is given by
$$
(H^2f)(t) = \frac{\sqrt{2\mu}}{i^{1/6}} \, e(\mu t^2) \intg f(x) e(2\mu tx) dx.
\tag3
$$ 
\endproclaim

\noindent
Therefore, $H^{-1} = H^5 = H^3 H^2$ and by (3) one gets the formula for the 
inverse hexic transform
$$
(H^{-1}f)(t) = \frac{\sqrt{2\mu}}{i^{1/6}} \, e(\mu t^2) \intg f(x) e(-2\mu tx) dx.
$$
One similarly gets a formula for the inverse cubic transform 
$(H^{-2}f) = (H^3Hf)(t) = (Hf)(-t)$.

\bigpagebreak

\remark{Remark 1}
It is interesting to note that although the Fourier transform has period four,
if it is composed with multiplication by a complex Gaussian, as in
(3), it can be made to have period three.  Though this may seem a little 
surprising, it can be shown that from the C*-algebra point of view it is 
not (see Section 3).
\endremark

\remark{Remark 2} 
By analogy with the fact that $e^{-\pi x^2}$ is invariant under the Fourier
transform, one can easily check that $e^{-\pi(\sqrt3-i)\mu x^2}$ is invariant under 
$H$ (and hence also under the cubic transform).  (It can be checked
that this is the only function among the Gaussian exponentials, up to scalars, that is 
invariant under the cubic or hexic transform.)  
The reason we referred to $\mu = 1/2$ as the ``ideal'' case is that in this case 
one has $\tfrac12(\sqrt3-i) = i^{-1/3}$ is of modulus 1 so that one has the invariant 
Gaussian $e^{-\pi i^{-1/3} x^2}$. 
\endremark

\remark{Remark 3}
The Fourier transform is a ``canonical'' transform in the sense that it intertwines the
translation and phase multiplication operators in the well known way.  Similarly,
the hexic transform, with $\mu=\tfrac12$, is also canonical.  In fact, letting 
$(T_xf)(t)=f(t-x)$ and $(E_xf)(t)=e(-xt)f(t)$, one checks the following relations 
(see \S3 below): 
$$
T_xH = HE_x, \qquad  E_xHT_x =  e(-\tfrac12x^2) T_xH.
$$
Since the group $\roman{SL}(2,\Bbb Z)$ is known to contain finite order elements only of 
orders 2, 3, 4, and 6, it follows that a periodic canonical transform can only have
these orders. This therefore gives us canonical transforms for each allowed order.
\endremark

\remark{Remark 4} 
It may be worthwhile investigating properties that the above transforms
have that are analogous to those that are well known to hold of the Fourier transform
(as for example in Rudin [\WR]).  For example, is there a multiplication $\#$ on the
space of Schwartz functions (or $L^1$ functions) such that $H(f\#g) = H(f)H(g)$?
(For the Fourier transform this multiplication is convolution.)  It may also be of
interest to explore the extension of the transforms $H$ and $H^2$ to $\Bbb R^n$, 
or even to locally compact Abelian groups.
\endremark

\medpagebreak

An application of Theorem 1 is the existence of finitely generated projective modules over
crossed products $6_\theta := A_\theta \rtimes_\rho \Bbb Z_6$ and 
$3_\theta := A_\theta \rtimes_{\rho^2} \Bbb Z_3$, where $A_\theta$ is the rotation 
C*-algebra and $\rho$ is the canonical order six automorphism on $A_\theta$ (see \S3).  
These modules will give rise to primary classes in the corresponding $K_0$-groups 
$K_0(6_\theta)$ and $K_0(3_\theta)$.
We write $6_\theta^\infty$ and $3_\theta^\infty$ for the respective canonical smooth 
dense *-subalgebras.  It is well known [\AC] that there are natural isomorphisms 
$K_*(6_\theta) = K_*(6_\theta^\infty)$ and $K_*(3_\theta) = K_*(3_\theta^\infty)$.
(See Section 3.)

\proclaim{Theorem 2}
Under the action (2), the Schwartz space $\sr$ is a finitely generated 
projective right module $\Cal M_6$ over $6_\theta^\infty$ (thus giving rise to a class in 
$K_0(6_\theta^\infty)$).
Similarly, under the action of the order three unitary given by (3),
the Schwartz space $\sr$ is a finitely generated projective right module $\Cal M_3$ over 
$3_\theta^\infty$ (thus giving rise to a class in $K_0(3_\theta^\infty)$).
Further, one has
$$
\vartau_*[\Cal M_6] \ = \ \frac\theta6, \qquad
\vartau_*[\Cal M_3] \ = \ \frac\theta3,
$$
for $j=0,1,\dots,5$, where $\vartau_*$ is the induced map by the canonical trace $\vartau$ on $K_0$.
\endproclaim

\bigpagebreak 

In [\BW], Buck and the author compute the Connes-Chern characters of the hexic and 
cubic modules $\Cal M_6, \Cal M_3$ and show that there are explicit injections 
$\Bbb Z^{10} \to K_0(6_\theta)$ and $\Bbb Z^8 \to K_0(3_\theta)$ for each $\theta>0$.
The author believes that, just as in the Fourier case [\SWb], these injections will 
turn out to be isomorphisms (at least for a dense $G_\delta$ set of $\theta$).  

The author wishes to thank George Elliott for making some helpful suggestions.

\bigpagebreak

\subhead \S2. PROOF OF THEOREM 1\endsubhead

We will make free use of the following identity 
$$
\intg\ e(Ax)\,e^{-\pi bx^2}\,dx \ = \
\frac1{\sqrt b} e^{-\pi A^2/b}
$$
which holds for $b,A\in \Bbb C,\ \text{Re}(b) > 0$, and $\sqrt b$ is the principal 
square root.

The theorem follows once we show that:

\medpagebreak

\itemitem{(A)} the set of Gaussians $f_{\alpha}(x) := e(-\alpha x) e^{-2\pi\mu x^2}$, where 
$\alpha\in\Bbb R$, is a total set in $L^2(\Bbb R)$,

\itemitem{(B)} $(H^3f_\alpha)(t) = f_\alpha(-t)$ for all $t,\alpha$, 

\itemitem{(C)} equality (3) holds for each $f_\alpha$.

\medpagebreak

\demo{Proof of {\rm(A)}} 
It is enough to show that if $g\in L^2(\Bbb R)$ is such that 
$$
\intg g(x) e(-\alpha x) e^{-2\pi\mu x^2} dx  = 0 
$$
for each $\alpha$, then $g=0$.  Setting $g(x) e^{-2\pi\mu x^2} = h(x)$ we note that 
$h$ is in $L^1(\Bbb R)$ since it is a product of two $L^2$ functions.  Hence one has
$$
0 = \intg h(x) e(-\alpha x) dx = \ft h(\alpha)
$$
for each $\alpha$.  Therefore, $\ft h = 0$ and hence $h=0$, i.e., $g=0$.  This proves (A)
and shows that the set of linear combinations of functions of the form $f_\alpha$ is
a dense subspace of $L^2(\Bbb R)$.
\enddemo

\bigpagebreak

\demo{Proof of \rm{(B)}} One has
$$
\align
(Hf_\alpha)(t) &= i^{1/6}\sqrt{2\mu} \intg  e(-\alpha x) e^{-2\pi\mu x^2} e(2\mu tx - \mu x^2) dx
\\
&= i^{1/6} \sqrt{2\mu} \intg e((2\mu t-\alpha)x) \, e^{-2\pi\mu(1+i)x^2} dx
\\
&=
\frac{i^{1/6}}{\sqrt{1+i}}  e^{-\pi(2\mu t-\alpha)^2/(2\mu(1+i))}.
\endalign
$$
Applying $H$ again gives
$$
\align
(H^2f_\alpha)(t) &= \frac{i^{1/3}\sqrt{2\mu}}{\sqrt{1+i}} 
\intg  e^{-\pi(2\mu x-\alpha)^2/(2\mu(1+i))} e(2\mu tx - \mu x^2) dx
\\
&=
\frac{i^{1/3}\sqrt{2\mu}}{\sqrt{1+i}} e(C)
\intg e((2\mu t+D)x) e^{-\pi\beta x^2} \, dx
\\
&=
\frac{i^{1/3}\sqrt{2\mu}}{\sqrt{(1+i)\beta}} e(C) e^{-\pi(2\mu t+D)^2/\beta}
\\
&=
\frac{i^{1/3}}{\sqrt i} e(C) e^{-\pi(2\mu t+D)^2/\beta}
\tag4
\endalign
$$
where
$$
\beta = \mu(1+i), \qquad C = \tfrac1{8\mu}(i+1)\alpha^2, \qquad
D = \tfrac{-(i+1)}2\alpha.
$$
A third iteration gives
$$
\align
(H^3f_\alpha)(t) &= \sqrt{2\mu} \ e(C) 
\intg e^{-\pi(2\mu x+D)^2/\beta} e(2\mu tx - \mu x^2) dx 
\\
&=
\sqrt{2\mu}\ e^{-\tfrac{\pi\alpha^2}{2\mu}} 
\intg e((2\mu t-i\alpha)x) e^{-2\pi\mu x^2} \, dx
\\
&=
e^{-\tfrac{\pi\alpha^2}{2\mu}} \,e^{-\pi(2\mu t-i\alpha)^2/(2\mu)}
\\
&=
e(\alpha t) \, e^{-2\pi\mu t^2}.
\endalign
$$ 
Therefore, $(H^3f_\alpha)(t) = f_\alpha(-t)$ is the usual flip map. Since this
holds for all $\alpha$, and $\{f_\alpha\}$ is a total set of functions in $L^2(\Bbb R)$, 
this relation holds for all $L^2$ functions on $\Bbb R$.   Hence $H^6$ is the identity.  
\enddemo

\bigpagebreak

\demo{Proof of \rm{(C)}}
The right hand side of (3) evaluated at $f_\alpha$ is
$$
\frac{\sqrt{2\mu}}{i^{1/6}} \, e(\mu t^2) \intg  e(-\alpha x) e^{-2\pi\mu x^2} e(2\mu tx) dx
= \frac1{i^{1/6}} \, e(\mu t^2) \ e^{-\pi(2\mu t-\alpha)^2/2\mu} 
$$
and it is easy to check that this is exactly (4), namely $(H^2f_\alpha)(t)$.  This completes
the proof of Theorem 1.  
\enddemo

\bigpagebreak

\subhead \S3. APPLICATION TO C*-ALGEBRAS \endsubhead

The following shows how by means of C*-algebras one can discover the above transforms.

Let $\theta>0,\, \lambda=e(\theta)$, and consider the rotation C*-algebra $A_\theta$ 
generated by unitaries $U,V$ satisfying $VU=\lambda UV$.  The (noncommutative) 
{\it hexic} transform of $A_\theta$ is the canonical order six automorphism 
$\rho$ defined by
$$
\rho(U)=V, \qquad \rho(V) = \lambda^{-1/2}U^{-1}V.
$$
Its square $\kappa := \rho^2$ is the canonical order three automorphism, which we call
the {\it cubic} transform, and $\rho^3$ is the usual flip automorphism studied in
great detail in [\BEEKa], [\BEEKb], and [\BK].
The corresponding crossed product $6_\theta := A_\theta \rtimes_\rho \Bbb Z_6$ is the
universal C*-algebra generated by unitaries $U,V,W$ enjoying the commutation relations
$$
VU=\lambda UV, \qquad WUW^{-1} = V, \qquad WVW^{-1} = \lambda^{-1/2}U^{-1}V, \qquad W^6 = I.
\tag5
$$
One may view the crossed product $3_\theta = A_\theta \rtimes_{\kappa} \Bbb Z_3$ as 
the C*-subalgebra of $6_\theta$ generated by $U,V$, and $W^2$.
We write $6_\theta^\infty$ and $3_\theta^\infty$ for their respective canonical smooth 
dense *-subalgebras. (For example, the elements of $6_\theta^\infty$ consist of sums of terms
of the form $aW^j$ where $a\in A_\theta^\infty$.)
Using Rieffel's Theorem 2.15 [\MRb] (with an appropriate lattice group in $\Bbb R \times 
\hat{\Bbb R}$) one obtains a smooth Heisenberg module structure on the Schwartz space $\sr$,
with $A_\theta^\infty$ acting on the right, given by
$$
(fU)(t) = f(t-\alpha), \qquad (fV)(t) = e(-\alpha t)f(t),
$$
where $\alpha = \sqrt\theta$.  
To extend this action so as to obtain a right $6_\theta^\infty$-module action on $\sr$, 
we need $W$ to act as an integral transform
$$
(fW)(t) = \intg f(x) K(x,t) dx,
$$
for suitable kernel function $K$, so that the relations (5) are satisfied.
Rewrite the commutation relations in (5) involving $W$ in the form
$$
WU = VW, \qquad VW = \lambda^{1/2}UWV, \qquad W^6 = I.
\tag6
$$
For the second of these relations one has
$$
(fUWV)(t) = e(-\alpha t) (fUW)(t) = e(-\alpha t) \intg f(x) K(x+\alpha,t) dx 
$$
and
$$
(fVW)(t) = \intg e(-\alpha x) f(x) K(x,t) dx.
$$
Hence, doing the same thing for the first relation in (6), one gets the relations
$$
K(x,t-\alpha) = e(-\alpha x)\, K(x,t), \qquad  
\lambda^{1/2}\, K(x+\alpha,t) = e(\alpha t-\alpha x)\, K(x,t).
$$
Now it is easy to check that the kernel function $K(x,t) = i^{1/6} e(tx - \tfrac12 x^2)$ 
(of the transform $H$ above with $\mu = \tfrac12$) satisfies these relations.
Therefore, by Theorem 1 we can define the right action of $W$ on $\sr$ by:
$$
(fW)(t) \ = \ i^{1/6} \intg f(x) e(tx-\tfrac12 x^2) dx,
$$
so that the three relations in (6) hold.
This, together with the above actions of $U,V$ gives rise to a right $6_\theta^\infty$ module
structure on $\sr$.  We shall denote this module by $\Cal M_6$ and call it the {\it hexic}
module.
Also by Theorem 1, one has the order three action
$$
(fW^2)(t) \ = \ i^{-1/6} e(\tfrac12 t^2) \intg f(x) e(-tx) dx 
\ = \ i^{-1/6} e(\tfrac12 t^2) \, \ft f(t),
$$
which makes $\sr$ into a right $3_\theta^\infty$-module---we call it the {\it cubic} module
and denote by $\Cal M_3$.
In view of Rieffel's inner product formulas [\MRb], the $A_\theta^\infty$-valued inner 
on the Heisenberg module $\sr$ can be defined by 
$$
\inner f g {A_\theta^\infty}  \ = \
\sum_{m,n} \ \inner f g {A_\theta^\infty} (m,n) \cdot V^n U^m,
$$
where $f,g \in \sr$ and 
$$
\inner f g {A_\theta^\infty} (m,n) \ = \ \int_{-\infty}^\infty
\overline{f(t+\alpha m)} g(t) \, e(-\alpha nt)\,dt.
$$
The associated $6_\theta^\infty$ and $3_\theta^\infty$-valued inner product 
$\inner {\ \ } {\ } {6_\theta^\infty}, \ \inner {\ \ } {\ } {3_\theta^\infty}$
are defined by symmetrization
$$
\inner f g {6_\theta^\infty} \ = \ 
\sum_{j=0}^5 \ \inner f {g W^{-j}} {A_\theta^\infty}\, W^j, \qquad
\inner f g {3_\theta^\infty} \ = \ 
\sum_{j=0}^2 \ \inner f {g W^{-2j}} {A_\theta^\infty}\, W^{2j}.
$$
With these inner products, and exactly as was done in the proof for the Fourier module
in [\SWa], one sees that the hexic and cubic modules are finitely generated projective, 
giving classes in the corresponding $K_0$-group, and therefore one obtains Theorem 2.

\bigpagebreak



\enddocument